\documentclass{amsart}
\usepackage{mathptmx}
\usepackage{newcent}
\usepackage{url}

\newif\ifdraft\draftfalse
\newif\ifcite\citefalse
\usepackage{amsfonts,amssymb, amsmath}
\usepackage{oldgerm}
\usepackage{amscd}

\usepackage[hyperindex]{hyperref}  % To Make ArXiv happy
\ifdraft\ifcite\usepackage{showkeys}\else\usepackage[notcite,notref]{showkeys}\fi\fi

\usepackage{xy}
\newdir{ >}{{}*!/-5pt/\dir{>}}

\newtheorem{theorem}[equation]{Theorem}

\newtheorem{lemma}[equation]{Lemma}
\newtheorem{conjecture}[equation]{Conjecture}
\newtheorem{corollary}[equation]{Corollary}

\theoremstyle{remark}
\newtheorem{definition}[equation]{Definition}
\theoremstyle{remark}

\newtheorem{remark}[equation]{Remark}

\numberwithin{equation}{section}

\def\define{\def}

\define\M{{\mathcal{M}}}
\define\V{{\mathcal{V}}}
\define\F{{\mathcal{F}}}
\define\W{{\mathcal{W}}}

\define\H{\mathcal{H}}

\define\C{\mathbb{C}}
\define\R{\mathbb{R}}
\define\Q{\mathbb{Q}}
\define\Z{\mathbb{Z}}
\define\HH{\mathbb{H}}
\define\lam{\Lambda^{-1,-1}}

\define\Gc{G_{\C}}

\define\Ext{\operatorname{Ext}}

\define\sing{\operatorname{sing}}
\define\ad{\rm ad\,}
\define\Ad{\rm Ad\,}
\define\pd{\partial}

\define\a{\alpha}
\define\b{\beta}
\define\g{\gamma}

\define\beq{\begin{equation}}
\define\eeq{\end{equation}}

\def\HH{\mathbb H}

\newcommand\dd{\mathcal D}

\newcommand\colim{\mathop{\mathrm{colim}}}

\newcommand\Aut{\mathop{\mathrm{Aut}}\nolimits}
\newcommand\GL{\mathop{\mathrm{GL}}\nolimits}
\newcommand\SL{\mathop{\mathrm{SL}}\nolimits}
\newcommand\lsl{\mathop{\mathfrak{sl}}\nolimits}
\newcommand\gl{\mathop{\mathfrak{gl}}\nolimits}

\newcommand\cl{\mathop{\mathrm{cl}}\nolimits}

\newcommand\MHS{\mathop{\mathrm{MHS}}\nolimits}

\newcommand\IH{\mathop{\mathrm{IH}}\nolimits}

\begin{document}
\bibliographystyle{plainurl}
\title{Zero loci of admissible normal functions with torsion singularities}
\author{Patrick Brosnan}
\address{Department of Mathematics\\
 The University of British Columbia\\
 1984 Mathematics Road}
\email{brosnan@math.ubc.ca}
\author{Gregory Pearlstein}
\address{Department of Mathematics\\ 
         Michigan State University\\ 
         East Lansing, MI 48824}
\email{gpearl@math.msu.edu}
%\date{March 20}
%\subjclass{}
% \keywords{}
%\dedicatory{}

\begin{abstract}
  We show that the zero locus of a normal function on a smooth complex
  algebraic variety S is algebraic provided that the normal function
  extends to a admissible normal function on a smooth compactification
  of S with torsion singularity.  This result generalizes our previous result 
for admissible normal functions on curves [arxiv:0604345]. It has also been obtained by
  M.~Saito using a different method in a recent preprint
  [arXiv:0803.2771v2].
\end{abstract}

\maketitle
% \tableofcontents

\section{Introduction} Let $H$ be a pure Hodge structure of weight $-1$
with integral structure $H_{\Z}$.  Then, the intermediate Jacobian of
$H$ is the complex torus
$
            J(H) = H_{\C}/(F^0 + H_{\Z})
$
where $F^*$ is the Hodge filtration of $H$.  If $\H$ is a variation of pure 
Hodge structure over a complex manifold $S$ with integral structure $\H_{\Z}$, 
the above construction produces a holomorphic bundle of complex tori 
$J(\H)\to S$ with fiber $J(\H)_s = J(\H_s)$ over $s$.  A normal function 
$\nu$ is a holomorphic section of $J(\H)$ which satisfies a version of 
Griffiths horizontality.  Therefore, as a holomorphic section of $J(\H)$, the 
locus of points $\mathcal Z$ where $\nu$ vanishes is a complex analytic 
subvariety of $S$.  Furthermore, we have the following conjecture
of Griffiths and Green:

\begin{conjecture}\label{conjecture:main}
Let $\nu$ be an admissible normal function~\cite{SaitoANF}
on a smooth complex algebraic variety $S$.  Then, the
zero locus $\mathcal Z$ of $\nu$ is an algebraic subvariety
of $S$.
\end{conjecture}

\par In analogy with the work of Cattani, Deligne and Kaplan~\cite{CDK}
on the algebraicity of the locus of a Hodge class, an unconditional 
proof of this conjecture provides evidence in support of the standard
conjectures on the existence of filtrations on Chow groups~\cite{GG2}.
In the case where $S$ is a curve, we gave an unconditional proof of 
\eqref{conjecture:main} in~\cite{BP}.  Other special cases in which 
\eqref{conjecture:main} is known are normal functions arising from cycles 
which are algebraically equivalent to zero and the case where the fibers 
of $J(\H)$ are Abelian varieties. % (representability of the Picard scheme)  
In this paper, we prove the following extension of~\cite{BP}:

\begin{theorem}\label{theorem:main-1} Let $\nu$ be an admissible 
normal function~\cite{SaitoANF} on a smooth complex algebraic 
variety $S$.  Assume that $S$ has a smooth compactification $\bar S$ 
such that $D=\bar S - S$ is a smooth divisor.  Then, the zero locus 
$\mathcal Z$ of $\nu$ is an algebraic subvariety of $S$.
\end{theorem}

\par The first step in the proof of Theorem \eqref{theorem:main-1} is to
replace $\nu$ be an admissible variation of mixed Hodge structure $\V$
with integral structure $\H_{\Z}$ and weight graded quotients $Gr^W_0$
and $Gr^W_{-1}=\H$.  This is possible by~\cite{SaitoANF}.  By a standard
construction of Deligne, the mixed Hodge structure on the fiber $\V_s$
defines a grading $Y(s)$ of the weight filtration of $\V_s$ which 
preserves the Hodge filtration.
The zero 
locus $\mathcal Z$ is then exactly the set of points where $Y(s)$ is defined 
over $\Z$.  
\vskip 3pt

\par In analogy with~\cite{BP}, the two key technical ingredients 
in the proof of Theorem \eqref{theorem:main-1} is the local normal form
of a variation of mixed Hodge structure along a normal crossing
divisor~\cite{P1} and the following lemma, which follows from
the the full strength of the 1-variable $\SL_2$-orbit theorem~\cite{P2}.

\begin{lemma}\label{lemma:main} Let $\Delta^r$ be a polydisk and
$D\subset\Delta^r$ be a smooth analytic hypersurface.  Let $\V$
be a variation of mixed Hodge structure over the complement of 
$D$ with weight graded quotients $Gr^W_0$ and $Gr^W_{-1}$.  
Assume that the monodromy $T=e^N$ of $\V$ about $D$ is unipotent.  
Then, for each point $p\in D$, the limit
$$
        \hat Y(p) = \lim_{s\to p}\, Y(s)
$$
exists, is contained in the kernel of $\ad N$ and has an
explicit description in terms of $N$ and the $\lsl_2$-splitting
of the limit mixed Hodge structure of $\V$ at $p$.
\end{lemma} 

\begin{remark} The limit mixed Hodge structure of $\V$ at $p$ depends upon
the choice of local coordinates of $\bar S$ at $p$.  However, because the
limit $\hat Y(p)$ belongs to the kernel of $\ad N$, it well defined
independent of the choice of local coordinates.
\end{remark}

\par Alternatively, instead of taking the limit of $Y(s)$ as $s$ accumulates
to $p\in D$ along a sequence of points in $S$, one can twist $Y(s)$ by 
$e^{-\frac{1}{2\pi i}\log(s)N}$ in analogy with the construction of the 
limit mixed Hodge structure.  This gives a corresponding grading 
$Y(p)$ which belongs to the kernel of $\ad N$ and has an explicit
description in terms of the limit mixed Hodge structure of $\V$ at $p$.
This is stated explicitly in Theorem $(4.15)$ of~\cite{P2}.

\par In terms of the grading $Y(s)$, the normal function $\nu$ is 
constructed as follows:  Let $Y_{\Z}$ be an integral grading of 
some reference fiber of $\V$.  Then, $Y_{\Z}$ extends to a 
multivalued, integral grading of the weight filtration of $\V$
over $S$.  Therefore, the difference $Y(s)-Y_{\Z}$ is a well
defined map from $\Z(0)$ into $J(\H_s)$ for each point
$s\in S$.  The normal function $\nu$ is the image of $1\in\Z(0)$
under this map.  This suggests setting
$$
     J(\H)_p = \Ext^1_{\MHS}(\Z(0),K) 
$$
where $K$ is the induced mixed Hodge structure on $\ker(N:H_{\C}\to H_{\C})$
and defining
$$
       \nu(p) = (Y(p) - Y_{\Z})(1)\in\bar J(\H)_p 
$$
where $T=e^N$ is the local monodromy of $\V$ at $p$ (assumed unipotent),
$F_{\infty}$ is the limit Hodge filtration of $\H$ at $p$, 
and $Y_{\Z}$ is an integral grading of the weight filtration
which is invariant under $T$.

\par In general, the existence of such a grading $Y_{\Z}$ is obstructed
by the class $\sigma_{\Z,p}(\nu)$ of $\nu$ in the finite group
\beq
       G=\frac{H_{\Z}\cap (T-1)(H_{\Q})}{(T-1)(H_{\Z})} \label{eq:finite-group}
\eeq  
In analogy with~\cite{GGK}, this allows one to construct a 
\lq\lq Neron model\rq\rq{} which graphs admissible normal functions on a 
neighborhood of $p$:  In general, the fibers of $\bar J(\H)$ can not patch
together to form a complex analytic space, since the dimension of 
$\bar J(\H)_p$ can be less than the dimension of $J(\H)_s$ for $s\in S$.
Nonetheless, $\bar J(\H)$ does carry a kind of generalized complex 
analytic structure (\lq\lq slit analytic space\rq\rq) which traces back 
to the  fundamental work of Kato and Usui compactification of period 
domains~\cite{KU}.  For recent work in this direction see~\cite{SaitoLetter} 
which uses the Neron model of~\cite{GGK} to give a proof of Theorem~\ref{theorem:main-1} independent of ours.  

\par Our original interest in the construction of the limits of normal
functions is rooted in the work of Griffiths and Green~\cite{GG}
on singularities of normal functions and the Hodge conjecture.  Very 
briefly, the idea of~\cite{GG} is to start with a smooth projective
variety $X$ of complex dimension $2n$ and a very ample line bundle
$L$ on $X$.  Let $|L| = \mathbb P H^0(X,L)$ and $S$ be the complement 
of the dual variety  $\hat X\subset |L|$ of $X$. Then (cf.~\cite{BFNP}), 
a primitive Deligne cohomology class $\zeta\in H^{2n}_{\mathcal D}(X,\Z(n))$ 
determines an admissible normal function $\nu$ on $S$ with cohomology class 
$\cl_{\Z}(\nu)\in H^1(S,\H_{\Z})$.  We then say that $\nu$ is singular on 
$|L|$ if there is a point $p\in\hat X$ such that
\beq
      \sigma_{\Z,p}(\nu) = \colim_{p\in U}\, \cl_{\Z}(\nu)|_{U\cap S}
                        \in\colim_{p\in U}\, H^1(S\cap U,\H_{\Z})
      \label{eq:singularity}
\eeq
is non-torsion, where the colimit is taken over all complex analytic
neighborhoods $U$ of $p$ in $|L|$.  The Hodge conjecture is then equivalent 
to the following statement~\cite{GG,BFNP}

\begin{conjecture} For each primitive, non-torsion Hodge class 
$\zeta\in H^{n,n}(X,\Z)$ there exists a positive integer $k$
such that $\nu$ is singular on $|L^k|$.
\end{conjecture}

\begin{remark} The definition of $\sigma_{\Z,p}(\nu)$ is valid for any
admissible normal function defined on the complement of a divisor 
$D\subset\bar S$.  The finite group \eqref{eq:finite-group} is exactly 
the torsion part of the cohomology group appearing in \eqref{eq:singularity}.  
In case where $D$ is a smooth divisor, admissibility forces 
$\sigma_{\Z,p}(\nu)$ to be torsion.  
\end{remark}

\par Simple examples show that, in general, unless $\sigma_{\Z,p}(\nu)=0$
the limit of $Y(s)$ along a holomorphic arc $\gamma$ through $p$ depends 
upon the multiplicities (assumed finite) of the intersection of $\gamma$ 
with the irreducible components of the (normal crossing) boundary divisor 
at $p$.  
However, we will show that if 
$\sigma_{\Z,p}=0$, the limit $Y(s)$ is independent of $\gamma$.
Furthermore, modulo one step which we shall defer to~\cite{BP2}, we obtain 
the following result:

\begin{theorem}\label{theorem:main-2} Let $\nu$ be an admissible normal 
function on a smooth complex algebraic variety $S\subset\bar S$.  Assume 
that $D=\bar S-S$ is a normal crossing divisor and that $\sigma_{\Z,p}(\nu)$ 
is torsion for every point $p\in D$.  Then, the zero locus of $\nu$ is an 
algebraic subvariety of $S$.
\end{theorem}

\begin{remark}
  In fact, the assumption that $D$ is a normal crossing divisor is not
  necessary.  To see this, suppose that we know the result in the case
  that $D$ is a normal crossing divisor.  Let $\nu$ be an admissible
  normal function on $\overline{S}$ which is smooth over $S$.  By
  Hironaka, we can find a resolution $p:\overline{T}\to \overline{S}$
  such that $p^{-1}S\to S$ is an isomorphism and
  $p^{-1}(\overline{S}\setminus S)$ is a normal crossing divisor.  It
  is easy to see that, if the singularity of $\nu$ is zero at every
  point in $\overline{S}$, then the singularity of the pullback of
  $\nu$ to $\overline{T}$ is zero on $\overline{T}$ as well.  Thus, by
  the theorem, the zero locus of $\nu$ on $S=\overline{T}\setminus
  p^{-1}S$ is algebraic.
\end{remark}

\begin{remark}  
  As mentioned above, Morihiko Saito has obtained an independent proof
  of Theorem~\ref{theorem:main-1}.  He also obtains
  Theorem~\ref{theorem:main-2}. ~See~\cite{SaitoLetter}.
\end{remark}

\section{Preliminary Results}

\subsection{Gradings and Splittings}

\par Let $V$ be a finite dimensional vector space over a field $k$ of 
characteristic zero, and 
$$
        0 = L_a \subseteq\cdots \subseteq L_i \subset L_{i+1}\subseteq\cdots
                \subseteq L_b =V
$$
be an increasing filtration of $V$ indexed by $\Z$.   Then, a grading
of $L$ is a semisimple endomorphism $Y$ of $V$ such that 
$$
        L_i = E_i(Y)\oplus L_{i-1}
$$
for each index $i$, where $E_i(Y)$ is the $i$-eigenspace of $Y$.  
Elements of $\GL(V)$ which preserve $L$ act on gradings of $L$ by 
the adjoint action:
$$
           g.Y = g Y g^{-1}
$$

\par Let $(F,W)$ be a mixed Hodge structure with Hodge filtration $F$
and weight filtration $W$.  Then~\cite{D1} % get a more exact reference
there exists a unique, functorial bigrading
$$
              V_{\C} = \bigoplus_{p,q}\, I^{p,q}
$$
of the underlying vector space $V_{\C}$ such that
\begin{itemize}
\item[{(a)}] $F^p = \oplus_{r\geq p}\, I^{r,s}$;        
\item[{(b)}] $W_k = \oplus_{r+s\leq k}\, I^{r,s}$;       
\item[{(c)}] $\bar I^{p,q} = I^{q,p} \mod \oplus_{r<q,s<p}\, I^{r,s}$.
\end{itemize}      
The associated Deligne grading $Y_{(F,W)}$ of $W$ is the semisimple
endomorphism of $V_{\C}$ which acts as multiplication by $p+q$ on
$I^{p,q}$. In particular, by properties $(a)$--$(c)$, if $g$ is an element 
of $\GL(V_{\R})$ which preserves $W$ then 
$$
      I^{p,q}_{(g.F,W)} = g.I^{p,q}_{(F,W)}
$$
with respect to the linear action of $\GL(V)$ on filtrations and subspaces.
Likewise, for $g$ as above
$
        Y_{(g.F,W)} = g.Y_{(F,W)}
$.

\par The mixed Hodge structure $(F,W)$ induces a mixed Hodge structure on 
the Lie algebra $\gl(V_{\C})$ with associated bigrading 
\beq
         \gl(V_{\C}) = \bigoplus_{p,q}\, \gl(V)^{p,q}     \label{eq:lie-bigrading}
\eeq  
Let $\lambda$ be an element of the subalgebra
$$
     \Lambda^{-1,-1} = \bigoplus_{a,b<0}\, \gl(V)^{a,b}
$$
Then, by properties $(a)$--$(c)$,
\beq
          I^{p,q}_{(e^{\lambda}.F,W)} = e^{\lambda}.I^{p,q}_{(F,W)}
                                                      \label{eq:equivariance}
\eeq
and hence $Y_{(e^{\lambda}.F,W)} = e^{\lambda}.Y_{(F,W)}$.

\begin{definition} A mixed Hodge structure $(F,W)$ is split over $\R$ if 
$\bar I^{p,q} = I^{q,p}$.
\end{definition}

\begin{lemma}\label{lemma:split-equiv} Let $(F,W)$ be a mixed Hodge structure.
Then, the following are equivalent:
\begin{itemize}
\item[{(a)}] $(F,W)$ is split over $\R$;
\item[{(b)}] $I^{p,q}_{(F,W)} = F^p\cap\bar F^q\cap W_{p+q}$;
\item[{(c)}] There exists a grading $Y$ of $W$ which preserves $F$ and
is defined over $\R$, in which case $Y=Y_{(F,W)}$.
\end{itemize}
\end{lemma} 

\par If $(F,W)$ is not split over $\R$ we can construct an associated split 
mixed Hodge structure $(e^{-i\delta}.F,W)$ as follows:

\begin{theorem} [Prop (2.20)~\cite{CKS}] There exists a unique real element
$\delta$ in $\Lambda^{-1,-1}$ such that
$  
           \bar Y_{(F,W)} = e^{-2i\delta}.Y_{(F,W)}
$.
Moreover, $\delta$ commutes with all $(r,r)$-morphism of 
$(F,W)$ and 
\beq
           (e^{-i\delta}.F,W)                      \label{eq:delta}
\eeq
is split over $\R$.
\end{theorem} 

\par Let $W$ be an increasing filtration of $V$ and $N$ be a nilpotent
endomorphism of $V$ which preserves $W$.  Assume that the relative
weight filtration~\cite{SZ} $M$ of $N$ and $W$ exists, and suppose
that there exists a grading $Y_M$ of $M$ which preserves $W$ and
satisfies the condition
\beq
           [Y_M,N] = -2N                             \label{eq:sl2-1}
\eeq
Let $Y$ be a grading of $W$ which preserves $M$, and 
$$
        N = N_0 + N_{-1} + \cdots
$$
be the decomposition of $N$ with respect to $\ad Y$ 
(i.e. $[Y,N_{-j}] = -jN_{-j}$). 

\begin{lemma}\label{lemma:Deligne-1} (Deligne~\cite{D2,KP}) 
Under the hypothesis of the previous paragraph, there exists a unique, 
functorial grading $Y=Y(N,Y_M)$ of $W$ which commutes with $Y_M$ such that:
\begin{itemize}
\item[{(a)}] $(N_0,H)$ is an $sl_2$-pair where $H=Y_M - Y$;
\item[{(b)}] If $(N_0,H,N_0^+)$ is the associated $sl_2$-triple
then $[N-N_0,N_0^+] = 0$.
\end{itemize}
\end{lemma}

\begin{corollary}\label{corollary:highest-weight} For $k>0$, $N_{-k}$ is 
either zero or a highest weight vector of weight $k-2$ with respect the 
representation of $sl_2$ constructed in the previous lemma.  In particular,
$N_{-1}=0$.
\end{corollary}

\subsection{Admissible nilpotent orbits} Let $\V\to S$ be a variation
of mixed Hodge structure over a complex manifold.  Then~\cite{P1,U}, 
in analogy with a variation of pure Hodge structure, a choice of reference
fiber $V$ for $\V$ allows us to represent $\V$ by a period map
$$
         \varphi:S\to\Gamma\backslash\M
$$
where $\M$ is a suitable classifying space of graded-polarized mixed Hodge
structures and $\Gamma$ is the image of the monodromy representation.
As in the pure case, the classifying space $\M$ is a submanifold of a 
suitable flag variety, and the period map $\varphi$ is holomorphic,
horizontal and locally liftable.  If $F:\tilde S\to\M$ is a lifting of 
$\varphi$ to the universal cover of $S$, then
$$
         \frac{\pd F^p}{\pd z_j}\subseteq F^{p-1},\qquad
         \frac{\pd F^p}{\pd \bar z_j}\subseteq F^p
$$
where $(z_1,\dots,z_r)$ are local holomorphic coordinates on $\tilde S$.

\par More precisely, let $Q_*$ be the graded-polarizations of $Gr^W$
and $\GL(V)^W$ denote the subgroup of $\GL(V)$ consisting of elements
which preserve $W$.  Define
$$
        G = \{\, g\in \GL(V)^W \mid Gr(g)\in \Aut_{\R}(Q_*)\,\}
$$
to be the subgroup of $\GL(V)^W$ consisting of elements which act by real
isometries of $Q$ on $Gr^W$.  Then, in analogy with the pure case,
$G$ acts transitively on $\M$ by biholomorphisms.  Likewise, we
have an embedding of $\M$ into its \lq\lq compact dual\rq\rq{}
$$
         \M = G/G^F \hookrightarrow \Gc/\Gc^F = \check\M
$$
where
$\Gc = \{\, g\in \GL(V)^W \mid Gr(g)\in \Aut_{\C}(Q_*)\,\}$, 
and $G^F$, $\Gc^F$ are the corresponding isotopy groups of some
point $F\in\M$.  The set of points $F\in\M$ for which the 
corresponding mixed Hodge structure $(F,W)$ is split over $\R$
is a homogeneous space for the Lie group $G_{\R} = G\cap \GL(V_{\R})$.
Define $\mathfrak g_{\R}$ and $\mathfrak g_{\C}$ to be the respective
Lie algebras of $G_{\R}$ and $\Gc$.

\par Let $\Delta\subset\C$ be the unit disk and $\V$ be a variation
of mixed Hodge structure on the complement $\Delta^*$ of the origin
with unipotent monodromy $T=e^N$.  Then, we have a commutative diagram
$$
\begin{CD}
      U @>F >> \M   \\
    @V s=e^{2\pi iz} VV @VVV \\
     \Delta^* @>\varphi >> \Gamma\backslash\M
\end{CD}
$$
Therefore, $\psi(z) = e^{-zN}.F(z):\Delta^*\to\check\M$ descends to a well
defined holomorphic map $\psi:\Delta^*$ into $\check\M$.  If 
$\V$ is admissible then ~\cite{SZ} 
\begin{itemize}
\item[{(a)}] $F_{\infty} = \lim_{s\to 0}\, \psi(s)\in\check\M$ exists;
\item[{(b)}] The relative weight filtration $M$ of $N$ and $W$ exists.
\end{itemize}
In this case~\cite{SZ},
\begin{itemize}
\item[{(i)}] $(F_{\infty},M)$ is a mixed Hodge structure relative to which 
$N$ is a $(-1,-1)$-morphism;
\item[{(ii)}] $(e^{zN}.F_{\infty},W)$ is an admissible nilpotent orbit.
\end{itemize}

\par For variations of mixed Hodge structure over a higher dimensional
base, Kashiwara defined admissibility via a curve test~\cite{Kashiwara}. 
In particular, if $\V$ is an admissible variation of mixed Hodge structure
defined on the complement of a normal crossing divisor with unipotent
monodromy transformations $T_j = e^{N_j}$, then the relative weight 
filtration $M(N_j,W)$ of $W$ and $N_j$ exists for each $j$.

\par The remainder of this section is devoted to the discussion of the 
1-variable $\SL_2$-orbit theorem~\cite{P2} which allows us to approximate
the nilpotent orbit $\theta(z)$ by an associated $\SL_2$-orbit $\hat\theta(z)$ 
arising from a representation $\rho:\SL_2(\R)\to G_{\R}$. We start by 
returning to Lemma~\eqref{lemma:Deligne-1}:

\begin{lemma}\label{lemma:Deligne-2} (Deligne~\cite{D2,KP}) Let 
$(F,N,W)$ define an admissible nilpotent orbit with relative weight 
filtration $M$.  Let $Y_M = Y_{(F,M)}$  and $Y=Y(N,Y_M)$ be the 
associated grading of Lemma \eqref{lemma:Deligne-1}.  Then, $Y$ preserves 
$F$.  If $(F,M)$ is split over $\R$ then $\bar Y = Y$.
\end{lemma}
\begin{proof} This follows from the functoriality of Deligne's grading
together with an explicit computation in the case where $(F,M)$ is 
split over $\R$.
\end{proof} 

\begin{definition} A mixed Hodge structure $(F,W)$ is of type $(I)$ if 
there exists an index $i$ such that $Gr^W_k = 0$ unless $k=i$, $i+1$.
\end{definition}

\begin{lemma} Every mixed Hodge structure of type $(I)$ is split over
$\R$.
\end{lemma}
\begin{proof} This follows directly from the short length of the 
weight filtration and property $(c)$ of Deligne's bigrading.
\end{proof}

\par Combining the above result, we now obtain a formula for 
$Y_{(e^{zN}.F,W)}$ along an admissible nilpotent orbit of type
$(I)$, i.e.  $(e^{zN}.F,W)$ is mixed Hodge structure of type $(I)$ 
for $Im(z)\gg 0$, when the associated limit mixed Hodge structure
$(F,M)$ is split over $\R$:

\begin{theorem}\label{corollary:Deligne-3} Let $(e^{zN}.F,W)$ be 
an admissible nilpotent orbit of type $(I)$.  Let $Y=Y(N,Y_M)$
be the associated grading of $W$ of Lemma \eqref{lemma:Deligne-2},
and suppose that $(F,M)$ is split over $\R$.  Then, for $Im(z)>0$:
\beq
          Y = Y_{(e^{zN}.F,W)}            \label{eq:split-orbit-formula}
\eeq
\end{theorem}
\begin{proof} The fact that $(e^{zN}.F,W)$ is a mixed Hodge structure for 
$Im(z)>0$ follows from the fact that $(F,M)$ is split over $\R$ and the
results of~\cite{CKS}.  By Corollary \eqref{corollary:highest-weight}
and the short length of $W$, $N_0 = N$.  Therefore, $Y$ preserves
$e^{zN}.F$ since $[Y,N]=0$ and $Y$ preserves $F$ by 
Lemma \eqref{lemma:Deligne-2}.  As $Y$ is defined over $\R$, 
\eqref{eq:split-orbit-formula} now follows from part $(c)$ of Lemma 
\eqref{lemma:split-equiv}.
\end{proof}

\define\gg{{\mathfrak g}}

\par The next result allows us to compute the asymptotic behavior of 
$Y_{(e^{zN}.F,W)}$ along an arbitrary admissible nilpotent orbit of 
type $(I)$.

\begin{theorem} ($\SL_2$-orbit theorem \cite{P2})\label{theorem:sl2-orbit}  
Let $(e^{zN}.F,W)$ be an admissible nilpotent orbit of type $(I)$ with 
relative weight filtration $M$.  Let $(\tilde F,M) = (e^{-i\delta}.F,M)$
denote Deligne's $\delta$-splitting \eqref{eq:delta} of $(F,M)$.  Then, 
there exists an element
$$
        \zeta\in\gg_{\R}\cap\ker(\ad N)\cap\lam_{(\tilde F,M)}
$$
and a distinguished real analytic function $\tilde g:(a,\infty)\to G_{\R}$
such that 
\begin{itemize}
\item[{(a)}] $e^{iyN}.F = \tilde g(y)e^{iyN}.\tilde F$ for $y>a$;       
\item[{(b)}] $\tilde g(y)$ and $\tilde g^{-1}(y)$ have convergent series 
expansions about $\infty$ of the form
\begin{eqnarray*}
       \tilde g(y) 
         &=& e^{\zeta}(1 + \tilde g_1 y^{-1} + \tilde g_2 y^{-2} + \cdots) \\
        \tilde g^{-1}(y) 
         &=& e^{-\zeta}(1+\tilde f_1 y^{-1} + \tilde f_2 y^{-2} + \cdots)
\end{eqnarray*}
with $\tilde g_k$, $\tilde f_k\in\ker (ad N)^{k+1}$;
\item[{(c)}] $\delta$, $\zeta$ and the coefficients $\tilde g_k$ are related by 
the formula
$$
         e^{i\delta} 
         = e^{\zeta}
           \left(1+\sum_{k>0}\, \frac{(-i)^k}{k!}(\ad N)^k \tilde g_k\right)
$$
\end{itemize}  
Let $(N_0,H,N_0^+)$ be the $sl_2$-triple determined by the $sl_2$-pair of 
Lemma \eqref{lemma:Deligne-2} and the nilpotent orbit $e^{zN}.\tilde F$.  
The constant $\zeta$ can be expressed as a universal Lie polynomial in the 
Hodge components $\delta^{r,s}$ of $\delta$ with respect to $(\tilde F,M)$.  
Likewise the coefficients $\tilde g_k$ and $\tilde f_k$ can be expressed 
as universal Lie polynomials in the Hodge components $\delta^{r,s}$ 
and $\ad N_0^+$.
\end{theorem}

\begin{remark} As noted in the proof of Corollary \eqref{corollary:Deligne-3},
for orbits of type $(I)$, $N=N_0$.
\end{remark}

\par For the purpose of computing the asymptotic behavior of the limit
grading in \S 3, it is useful to renormalize the $\SL_2$-orbit theorem
as follows: Let 
$$
      g(y) = \tilde g(y)e^{-\zeta},\qquad \hat F = e^{\zeta}.\tilde F
$$
Then, $e^{iyN}.F = g(y)e^{iyN}.\hat F$ since $[N,\zeta] =0$.  The mixed 
Hodge structure $(\hat F,M)$ is split over $\R$ since $(\tilde F,M)$
is split over $\R$ and $\zeta\in\mathfrak g_{\R}$.  Moreover,
$$
            e^{-\xi} = e^{\zeta}e^{-i\delta}\in\exp(\lam_{(\hat F_{\infty},M)})
$$
commutes with $N$ and is a universal polynomial in the Hodge components of 
$\delta$.  Likewise, the coefficients
$$
           g_k = \Ad(e^{\zeta})\tilde g_k 
$$
of the series expansion 
$$
       g(y) = 1 + \sum_{k>0} g_k y^{-k}
$$
are universal polynomials in the Hodge components of $\delta$ and 
$\ad N_0^+$, and satisfy the identity $g_k\in\ker(\ad N_0)^{k+1}$. 

\begin{definition} Let $(e^{zN}.F,W)$ be a nilpotent orbit of type
$(I)$.  Then,
$$
        \hat F = e^{-\xi}.F
$$
is the $sl_2$-splitting of $(F,M)$.
\end{definition}

\begin{remark} By virtue of the fact that $\zeta$ is given by a universal
polynomial in the Hodge components of $\delta$, the $sl_2$-splitting is
defined for any mixed Hodge structure.  The formula is as follows~\cite{KNU}:
Write the Campbell--Baker--Hausdorff formula as
$e^{\alpha}e^{\beta} = e^{H(\alpha,\beta)}$.  Then, $\delta$ and $\xi$ are
related by the formula
$$
         \delta = H(\xi,-\bar\xi)/2\sqrt{-1}
$$
\end{remark}
 
\subsection{Local normal form}  Let $\Delta^r$ be a polydisk with local
coordinates $(s_1,\dots,s_r)$ and $\V$ be an admissible variation of mixed 
Hodge structure on the complement of the divisor $s_1\cdots s_r=0$ with
unipotent monodromy $T_j = e^{N_j}$ about $s_j=0$.  Then, the $I^{p,q}$'s 
of the limit mixed Hodge structure $(F_{\infty},M)$ define a vector space 
complement
$$
      {\mathfrak q} = \bigoplus_{a<0} \gg^{a,b}
$$
to the isotopy algebra $\gg_{\C}^{F_{\infty}}$, and hence by admissibility,
near $p$ we can write the Hodge filtration of $\V$ as
\beq
      F(z) = e^{\sum_j\, z_j N_j}e^{\Gamma(s)}.F_{\infty}       \label{eq:2.1}
\eeq
where $\Gamma(s)$ is a $\mathfrak q$-valued function which vanishes at $s=0$ 
and $s_j = e^{2\pi i z_j}$.  By horizontality,
\beq
      \frac{\pd}{\pd z_j} F^p(z) \subseteq F^{p-1}(z)      \label{eq:2.2}
\eeq       

\par Let $\wp_a =\oplus_b\,\gg^{a,b}$ and note that:
\begin{itemize}
\item[{(i)}] $\mathfrak q = \oplus_{a<0}\,\wp_a$;
\item[{(ii)}] $[\wp_a,\wp_b] \subseteq \wp_{a+b}$;
\item[{(iii)}] $N_j\in\wp_{-1}$.
\end{itemize}
Inserting \eqref{eq:2.1} into \eqref{eq:2.2} it then follows 
that
\beq
      \Ad(e^{-\Gamma(s)})N_j 
        + 2\pi i s_j e^{-\Gamma(s)}\frac{\pd}{\pd s_j} e^{\Gamma(s)}
      \in\wp_{-1}                            \label{eq:2.3}
\eeq
Taking the limit as $s_j\to 0$ in \eqref{eq:2.3} 
it then follows by (i)--(iii) that
\beq
       [\Gamma^{(j)},N_j] = 0                              \label{eq:2.4}
\eeq
where $\Gamma^{(j)}$ denotes the restriction of $\Gamma(s)$ to the 
slice $s_j = 0$.

\begin{remark} In the pure case, this result is due to Cattani
and Kaplan~\cite{CK}. % Cattani and Kaplan's Luminy paper
\end{remark}

\begin{remark} The results of this section remains valid in the
case where $\V$ is a variation over $\Delta^{*a}\times\Delta^b$
upon setting $N_j =0$ for $j=a+1,\dots,b$.
\end{remark}

\subsection{Intersection Cohomology}

\par Let $\mathcal A_{\Q}$ be a local system of $\Q$-vector spaces over
a product of punctured disks $\Delta^{*r}$ with unipotent monodromy.
Let $A_{\Q}$ be a reference fiber of $\mathcal A_{\Q}$ and
$
           N_j\in{\rm Hom}(A_{\Q},A_{\Q})
$
denote  the monodromy logarithm of $\mathcal A_{\Q}$ about the j'th disk.  
Then, because the $N_j$'s commute, the vector spaces
\beq
         B^p(N_1,\dots,N_r;A_{\Q})
            = \bigoplus_{1\leq j_1 < \cdots < j_p \leq r}
              N_{j_1} N_{j_2}\cdots N_{j_p}(A_{\Q})       \label{eq:B-cmplx}
\eeq
form a complex with respect to the differential $d$ which acts on
the summands of \eqref{eq:B-cmplx} by the rule
\beq
             d:N_{j_1}\cdots \hat{N}_{j_q}\cdots N_{j_p}(A_{\Q})
             \stackrel{(-1)^{q-1}N_{j_q}}{\longrightarrow}
               N_{j_1}\cdots N_{j_p}(A_{\Q})
                                                        \label{eq:B-diff}
\eeq
Let $j:\Delta^{*r}\to\Delta^r$ be a holomorphic embedding (the open inclusion) 
of $\Delta^{*r}$ in a product of disks $\Delta^r$ and define
$$
        \IH^p(\Delta^r,\mathcal A_{\Q}) = \HH^p(\Delta^r,j_{!*}\mathcal A_{\Q})
$$
Then, by~\cite{D2, GGM} or \cite{KK}[Corollary 3.4.4]:
$H^p(B^*(N_1,\dots,N_r;A_{\Q})\cong \IH^p(\Delta^r,\mathcal A_{\Q})$.

The following result follows from Theorem~\cite{BFNP}[Lemma 2.1.8].
Here we give a proof that is more in the spirit of the calculations
done in this paper.

\begin{theorem}\label{theorem:ih-exact} Let $\V\to\Delta^{*r}$ be an
admissible variation of graded-polarizable mixed Hodge structure with
unipotent monodromy which is an extension of $\Q(0)$
by a variation of Hodge structure $\H$ of pure weight $-1$.  Then,
the associated short exact sequence
\beq
      0 \to \H_{\Q}\stackrel{\a}{\to} \V_{\Q}
            \stackrel{\b}{\to} \Q(0) \to 0               \label{eq:ih1}
\eeq
induces a long exact sequence
$$
       \cdots\to\IH^{p-1}(\Delta^r,\Q(0))\stackrel{\pd}{\to}
       \IH^p(\Delta^r,\H_{\Q})
       \stackrel{\a_*}{\to} \IH^p(\Delta^r,\V_{\Q})
       \stackrel{\b_*}{\to} \IH^p(\Delta^r,\Q(0))
       \to\cdots
$$
in intersection cohomology.
\end{theorem}
\begin{proof} Let $\V_{\Q}$ underlie an admissible extension of $\Q(0)$ by a 
variation of pure Hodge structure $\mathcal H$ of weight $-1$ and
$B^*(H_{\Q})$, $B^*(V_{\Q})$, $B^*(\Q(0))$
denote the associated complexes \eqref{eq:B-cmplx} defined by the
local monodromy.  By abuse of notation, let $\IH^p(H_{\Q})$, etc.
denote the cohomology of the the corresponding complex.
In particular, since each $N_j$ acts
trivially on $\Q(0)$, 
$$
      \IH^0(\Q(0)) = \Q(0),\qquad \IH^p(\Q(0)) = 0,\quad p>0
$$
Furthermore, since $N_j$ acts trivially on $\Q(0)$ and
$Gr^W_0(V_{\Q}) \cong \Q(0)$ it then follows that
$$
       N_j(V_{\Q}) \subset W_{-1}(V_{\Q}) \cong H_{\Q}
$$
By the existence of the relative weight filtration $M_j = M(N_j,W)$
and the short length of $W$ it then follows\cite{SZ} that
$$
        N_j(V_{\Q}) = N_j(H_{\Q})
$$
and hence $B^p(V_{\Q}) = B^p(H_{\Q})$ for $p>0$.
Consequently,
$$
        \IH^p(V_{\Q}) = \IH^p(H_{\Q}),\qquad p>1
$$
Combining the above results, we therefore obtain the
exactness of 
$$
\cdots\to \IH^{p-1}(\Q(0))\stackrel{\pd}{\to}
       \IH^p(\H_{\Q})
       \stackrel{\a_*}{\to} \IH^p(\V_{\Q})
       \stackrel{\b_*}{\to} \IH^p(\Q(0))
       \to\cdots
$$
for $p>1$.

\par Thus, in order to complete the proof, it remains to prove the exactness 
of the sequence
\beq
     0\to\IH^0(H_{\Q})\to\IH^0(V_{\Q})\to\IH^0(\Q(0))
       \stackrel{\pd}{\to}\IH^1(H_{\Q})\to\IH^1(V_{\Q})\to 0
                                       \label{eq:long-exact}
\eeq
By definition,
$$
      \IH^0(H_{\Q}) = \cap_j\ker(\left.N_j\right|_{H_{\Q}}),\qquad
      \IH^0(V_{\Q}) = \cap_j\ker(N_j)
$$
and hence the map
$
      \IH^0(H_{\Q}) \to \IH^0(V_{\Q})
$
is injective.

\par To see that \eqref{eq:long-exact} is exact at
$\IH^0(V_{\Q})$ observe that since $H_{\Q} = W_{-1}(V_{\Q})$
and $\Q(0) = Gr^W_0(V_{\Q})$, the image of $\IH^0(H_{\Q})$ in
$\IH^0(V_{\Q})$ is exactly the kernel of the map
\beq
             Gr^W_0:\IH^0(V_{\Q}) \to \IH^0(\Q(0))         \label{eq:ih-proj}
\eeq

\par For any class $[v]\in IH^0(\Q(0))$,
\beq
        \pd[v] = (N_1(v),\dots,N_r(v)) \mod dB^0(H_{\Q})    \label{eq:connect}
\eeq
where $v\in V_{\Q}$ is any element which projects onto
$[v]\in\Q(0) =  Gr^W_0(V_{\Q})$.  In particular, $\pd[v] = 0$
if and only if there exists $h\in H_{\Q} = W_{-1}(V_{\Q})$
such that
$$
             N_j(v) = N_j(h)
$$
for all $j$.  In this case,
$
            v_o = v - h
$
defines an element of $\IH^0(V_{\Q})$ which projects onto 
of $[v]\in \IH^0(\Q(0))$ under \eqref{eq:ih-proj}.  As such,
\eqref{eq:long-exact} is exact at $\IH^0(\Q(0))$.

\par To see that \eqref{eq:long-exact} is exact at $IH^1(H_{\Q})$
suppose that
$$
           (N_1(h_1),\dots,N_r(h_r)),\qquad h_j\in H_{\Q}
$$
represents a class $\eta\in\IH^1(H_{\Q})$ which maps to zero under
inclusion in $IH^1(V_{\Q})$.  Then, there exists a vector $v\in V_{\Q}$
such that
$$
            N_j(h_j) = N_j(v)
$$
for all $j$,  If $v\in H_{\Q}$ then $\eta = 0$.  Otherwise, $[v]$
defines a non-zero class in $\IH^0(\Q(0))$ such that
$\eta = \pd[v]$. Finally, to verify the surjectivity of the map
$$
       \IH^1(H_{\Q})\to \IH^1(V_{\Q})
$$
note that $B^p(H_{\Q}) = B^p(V_{\Q})$ for $p>0$ and
$dB^0(H_{\Q})\subseteq dB^0(V_{\Q})$.
\end{proof}

\begin{definition} Let $[1]$ be the class of $1$ in $\IH^0(\Q(0))$ and 
$\pd:\IH^0(\Q(0))\to\IH^1(\H_{\Q})$ be the connecting homomorphism.
Then, $\sing_p(\nu) = \pd 1$.
\end{definition}

\begin{remark} The results of this section remain valid upon replacing
$\Q$ by $\R$.
\end{remark}

\subsection{Invariant Grading} Let $\nu$ be an admissible normal function
over a product of punctured disks $\Delta^{*r}\subset\Delta^r$ with 
associated variation of mixed Hodge structure $\V$, reference fiber $V$ 
and nilpotent orbit $\theta({\bf z})=e^{\sum_j\, z_j N_j}.F_{\infty}$.  Let 
$H_{\R} = Gr^W_{-1} V_{\R}$ and $0=(0,\dots,0)\in\Delta^r$.  Let 
$\hat\theta({\bf z}) = e^{\sum_j\, z_j N_j}.\hat F$ be the split 
($sl_2$ or Deligne's $\delta$) orbit attached to  $\theta({\bf z})$. 
Let $\hat Y_M$ denote the corresponding grading of $(\hat F,M)$ where $M$ 
is the relative weight filtration of $W$ and the monodromy cone 
$$
      \mathcal C = \{\, \sum_j\, a_j N_j \mid a_j >0 \,\}
$$
Let $\hat Y = Y_{(e^{iN}.\hat F,W)}$, where $N= \sum_j\, N_j$  Then, by 
Lemma \eqref{lemma:Deligne-2}, $\hat Y$ is real, preserves $\hat F$ and 
commutes with $N$

\par Suppose that $\sing_0(\nu)=0$ and let $e_o$
be the element of $E_0(\hat Y)$ which projects to $1\in\R(0)$.
Define $e_j = N_j(e_0)$.  Then, by \eqref{eq:connect}
$$
       (e_1,\dots,e_r)\in B^1(H_{\R})
$$
is a representative of $\sing_p(\nu)$. Therefore, since $\sing_p(\nu) = 0$ 
there is an element $f\in H_{\R} =B^0(H_{\R})$ such that $e_j = N_j(f)$.  
Furthermore, since $e_j = N_j(e_0)$ and $e_0\in \hat F^0$ we have
$e_j \in \hat F^{-1}$.  Therefore, by strictness of morphisms of MHS, 
we can assume $f\in\hat F^0$.  Then, 
$$
       e_0-f = e^{iN}.(e_0-f) \in e^{iN}.\hat F^0
$$
Consequently, $e_0-f$ belongs to $I^{0,0}_{(e^{iN}.\hat F,W)}$ since
$e_0-f$ is real.  On the other hand, by theorem of Deligne
$e_0$ belongs to $I^{0,0}_{(e^{iN}.\hat F,W)}$.  Since $Gr^W_0$
has rank 1, it then follows that $f=0$.

\begin{corollary} $e_0\in\ker(N_j)$ for all $j$.
\end{corollary}

\begin{corollary}\label{corollary:inv-grading-split} If $\sing_0(\nu)=0$
then 
$$
         Y_{(e^{\sum_j\, z_j N_j}.\hat F,W)} = Y_{(e^{iN}.\hat F,W)}
$$
for $Im(z_1),\dots, Im(z_r)>0$.
\end{corollary}
\begin{proof} Since $e_0\in\ker(N_j)$ for all $j$, the grading 
$\hat Y = Y_{(e^{iN}.\hat F,W)}$ commutes with $N_1,\dots,N_r$.
Therefore, $\hat Y$ is real and preserves $e^{\sum_j\, z_j N_j}.\hat F$
since $\hat Y$ preserves $\hat F$, and hence is the Deligne grading
of $(e^{\sum_j\, z_j N_j}.\hat F,W)$.
\end{proof}

\par Let $Y_M = Y_{(F,M)}$.  Then, $Y_M = e^{i\delta}.\hat Y_M$ and 
$Y = Y(N,Y_M) = e^{i\delta}.\hat Y$.  Therefore, since $[\delta,N_j] = 0$
for all $j$, we have:

\begin{corollary}\label{corollary:inv-grading} $[Y,N_j] = 0$ for all $j$.
\end{corollary}

\begin{definition}\label{definition:inv-grading} If $\sing_0(\nu)=0$ we define 
$Y_{\infty}$ to be the grading $Y$ of Corollary \eqref{corollary:inv-grading}.  
In particular, since $Y_{\infty}$ commutes with $N_1,\dots,N_r$, it is 
independent of the choice of local coordinates used in its construction.
\end{definition}

\section{Limit Gradings} Let $D\subset\bar S$ be a smooth divisor, $p\in D$
and $\Delta^r$ be an analytic polydisk in $\bar S$
containing $p$.  Pick local coordinates $(s_1,\dots,s_r)$ on $\Delta^r$
such that $D\cap\Delta^r$ is given by $s_1 = 0$.  Represent $\nu$ by
an admissible variation of mixed Hodge structure $\V$ over 
$\Delta^*\times\Delta^{r-1}$ with weight graded quotients
$Gr^W_0 =\Z(0)$ and $Gr^W_{-1} = \mathcal H$.  Assume that the local monodromy 
of $\V$ about $D$ is given by a unipotent transformation $T=e^N$.  Let
$$
       F(z;s_2,\dots,s_r):U\times\Delta^{r-1}\to\M
$$
be a lifting of the local period map of $\V$ where $U$ is the upper 
half-plane.

\par Let 
$$
          F(z;s_2,\dots,s_r) = e^{zN}e^{\Gamma(s)}.F_{\infty}
$$
be the local normal form of the period map of $\V$ at $p$.
Let $\Gamma_0(s) = \Gamma(0,s_2,\dots,s_r)$ and 
$$
        F_{\infty}(s_2,\dots,s_r) = e^{\Gamma_0(s)}.F_{\infty}
$$
Let $W$ be the weight filtration of $\V$, $M$ be the relative
weight filtration of $N$ and $W$.  Then, 
$$
      \theta(z;s_2,\dots,s_r) = e^{zN}.F_{\infty}(s_2,\dots,s_r)
$$
is an admissible nilpotent orbit in 1-variable which depends 
complex analytically upon the parameters $(s_2,\dots,s_r)\in\Delta^{r-1}$.
Let
$$
          (\hat F_{\infty}(s_2,\dots,s_r),M) 
          = (e^{-\xi(s_1,\dots,s_r)}.F_{\infty}(s_2,\dots,s_r),M)
$$
denote the $sl_2$-splitting of $(F_{\infty}(s_2,\dots,s_r),M)$.
Then, $\xi$ is real analytic in $(s_2,\dots,s_r)$ since it
is given by universal Lie polynomials in the Hodge components
of Deligne's $\delta$-splitting of $(F_{\infty}(s_2,\dots,s_r),M)$.

\par By the $\SL_2$-orbit theorem \eqref{theorem:sl2-orbit}
$$
       \theta(iy;s_2,\dots,s_r)
       = g(y;s_2,\dots,s_r)
         e^{iy N}.\hat F_{\infty}(s_2,\dots,s_r)
$$
where 
$$
      g(y;s_2,\dots,s_r) 
      = (1 + \sum_{k>0}\, g_k(s_2,\dots,s_r) y^{-k})
$$ 
belongs to $G_{\R}$ and the coefficients $g_k(s_2,\dots,s_r)$ 
are real analytic in $(s_2,\dots,s_r)$ since they are given by universal 
Lie polynomials.  

\par We now derive an asymptotic formula for $Y_{(F(z;s_2,\dots,s_r),W)}$.
Write $z = x+iy$.  Then, 
\begin{eqnarray*}
     Y_{(F(z;s_2,\dots,s_r),W)}
     &=& Y_{(e^{xN}e^{iyN}e^{\Gamma(s)}.F_{\infty},W)}     \\
     &=& e^{xN}.Y_{(e^{iyN}e^{\Gamma(s)}e^{-\Gamma_0(s)}e^{\Gamma_0(s)}.
        F_{\infty},W)}                                       \\
     &=& e^{xN}.Y_{(e^{iyN}e^{\Gamma(s)}e^{-\Gamma_0(s)}.
        F_{\infty}(s_2,\dots,s_r),W)}                        
\end{eqnarray*}
Let $e^{\Gamma_1(s)} = e^{\Gamma(s)}e^{-\Gamma_0(s)}$ and note that 
$s_1 | \Gamma_1$ in $\mathcal O(\Delta^r)$.  Then, 
\begin{eqnarray*}
     Y_{(F(z;s_2,\dots,s_r),W)}
     &=& e^{xN}.Y_{(e^{iyN}e^{\Gamma_1(s)}.F_{\infty}(s_2,\dots,s_r),W)} \\
     &=& e^{xN}.Y_{(\Ad(e^{iyN})(e^{\Gamma_1(s)}).\theta(iy;s_2,\dots,s_r),W)}  \\
     &=& e^{xN}.Y_{(\Ad(e^{iyN})(e^{\Gamma_1(s)})
          g(y;s_2,\dots,s_r)e^{iy N}.\hat F_{\infty}(s_2,\dots,s_r),W)} 
\end{eqnarray*}

\par Let $ F_o(s_2,\dots,s_r) = e^{iN}.\hat F_{\infty}(s_2,\dots,s_r)$ and
$$
         Y_1 = Y_{(F_o(s_2,\dots,s_r),W)}, \qquad
         Y_2 = Y_{(\hat F_{\infty}(s_2,\dots,s_r),M)}
$$
Then, by Corollary \eqref{corollary:Deligne-3}:
$$
      Y_{(e^{iyN}.\hat F_{\infty}(s_2,\dots,s_r),W)} = Y_1
$$
Likewise, since $N$ and $H = Y_2 - Y_1$ is an $sl_2$-pair:
$$
      e^{iyN}.\hat F_{\infty}(s_2,\dots,s_r)
      = y^{-H/2}.F_o(s_2,\dots,s_r)
$$
Note that $Y_1$ and $H$ depend real analytically on $(s_2,\dots,s_r)$.        

\begin{lemma} Let $\gamma(y) =\Ad(e^{-iyN})g(y;s_2,\dots,s_r)$.  Then,
$\lim_{y\to\infty}\, \g(y)$ exists, and is real analytic in $(s_2,\dots,s_r)$.
\end{lemma}
\begin{proof} This follows directly from the fact that 
$g_k(s_2,\dots,s_r)$ is real analytic in $(s_2,\dots,s_r)$ 
and $g_k\in\ker(\ad N)^{k+1}$. 
\end{proof}

\par Returning to the calculation of $Y_{(F(z;s_2,\dots,s_r),W)}$, 
and abbreviating $g(y;s_2,\dots,s_r)$ to $g(y)$, we have
$$
\aligned
     Y_{(F(z;s_2,\dots,s_r),W)}
     &= e^{xN}.Y_{(\Ad(e^{iyN})(e^{\Gamma_1(s)})
          g(y)e^{iy N}.\hat F_{\infty}(s_2,\dots,s_r),W)}         \\ 
     &= e^{xN}.Y_{(g(y)e^{iyN}\g^{-1}(y)e^{\Gamma_1(s)}
          \g(y).\hat F_{\infty}(s_2,\dots,s_r),W)}    
\endaligned
$$
Let $e^{\Gamma_2} = \Ad(\g^{-1}(y))e^{\Gamma_1}$ and recall that 
$s_1 | \Gamma_1$.  Therefore,  
$$
\aligned
     Y_{(F(z;s_2,\dots,s_r),W)}
     &= e^{xN}.Y_{(g(y)e^{iyN}e^{\Gamma_2(s)}.
                       \hat F_{\infty}(s_2,\dots,s_r),W)}   \\
     &= e^{xN}g(y)y^{-H/2}.Y_{(e^{iN}\Ad(y^{H/2})(e^{\Gamma_2}).
         \hat F_{\infty}(s_2,\dots,s_r),W)}    
\endaligned
$$
where $\Ad(y^{H/2})\Gamma_2$ can be uniformly bounded by a constant
times $y^c e^{-2\pi y}$ as $y\to\infty$ for some constant $c$.  
\vskip 3pt

\par We now prove Lemma \eqref{lemma:main} of the introduction, which
we shall use in the next section to prove the algebraicity of the zero
locus $\mathcal Z$. Modulo our discussion of dependence on parameters,
this essentially the same calculation use to prove the existence of
the limit grading in~\cite{BP}.

\begin{theorem}\label{theorem:limit-1} Let $(s_1(m),\dots,s_r(m))$ be a 
sequence of points in $\Delta^*\times\Delta^{r-1}$ which converges to 
$(0,s_2,\dots,s_r)$ as $m\to\infty$.  Let $(z(m),s_2(m),\dots,s_r(m))$
be a lifting of this sequence to $U\times\Delta^{r-1}$ with the
real part of $z$ restricted to an interval of finite length.  
Then,
$$
     \lim_{m\to\infty}\, Y_{(F(z(m);s_2(m),\dots,s_r(m)),W)}
      = Y_{(e^{iN}.\hat F_{\infty}(s_2,\dots,s_r),W)}.
$$
\end{theorem}
\begin{proof} Suppress the dependence of $(z(m);s_2(m),\dots,s_r(m))$ on 
$m$.  By the previous results:
\beq
     Y_{(F(z;s_2,\dots,s_r),W)}
     = e^{xN}g(y)y^{-H/2}.Y_{(e^{iN}\Ad(y^{H/2})(e^{\Gamma_2(s)}).
         \hat F_{\infty}(s_2,\dots,s_r),W)}  \label{eq:untwisted-limit-1}
\eeq
where $\Ad(y^{H/2})(e^{\Gamma_2(s)}$ is uniformly bounded by some
constant times $y^c e^{-2\pi y}$.  Therefore,
\beq
         Y_{e^{iN}\Ad(y^{H/2})(e^{\Gamma_2(s)}).
            \hat F_{\infty}(s_2,\dots,s_r),W)}
         = Y_{(e^{iN}.\hat F_{\infty}(s_2,\dots,s_r),W)} + \a 
            \label{eq:untwisted-limit-2}     
\eeq
where $\a$ is uniformly bounded by $y^c e^{-2\pi y}$.  The result
now follows by inserting \eqref{eq:untwisted-limit-2} into 
\eqref{eq:untwisted-limit-1} and taking the limit at $m\to\infty$,
since $H=H(s_2,\dots,s_r)$ commutes with $Y_1(s_2,\dots,s_m)$.
\end{proof}

\par In order to construct the limit normal function we need the following
analog of Theorem \eqref{theorem:limit-1} where we twist the grading
$Y_{(\F,\W)}$ by $e^{-\frac{1}{2\pi i}\log(s_1) N}$.  Again, modulo dependence
on parameters, this is really just a glorified version of Theorem
$(4.15)$ in~\cite{P2}.

\begin{theorem}\label{theorem:limit-2} Let $(s_1(m),\dots,s_r(m))$ be a 
sequence of points in $\Delta^*\times\Delta^{r-1}$ which converges to 
$(0,s_2,\dots,s_r)$ as $m\to\infty$.  Let $(z(m),s_2(m),\dots,s_r(m))$
be a lifting of this sequence to $U\times\Delta^{r-1}$ with the
real part of $z$ restricted to an interval of finite length.  
Then,
$$
     \lim_{m\to\infty}\, e^{-zN}.Y_{(F(z(m);s_2(m),\dots,s_r(m)),W)}
      = Y(N,Y_{(F_{\infty}(s_2,\dots,s_r),M)}) 
$$
where  $Y(N,Y_{(F_{\infty}(s_2,\dots,s_r),M)})$ is the grading of 
Lemma \eqref{lemma:Deligne-2}.
\end{theorem}
\begin{proof} We repeat the argument of the proof of 
Theorem \eqref{theorem:limit-1} to obtain
\begin{eqnarray}
        e^{-zN}.Y_{(F(z;s_2,\dots,s_r),W)}
        &=& e^{-iyN}g(y)y^{-H/2}.
            (Y_{(e^{iN}.\hat F_{\infty}(s_2,\dots,s_r),W)} + \a)    \nonumber \\
        &=& e^{-iyN}\tilde g(y)e^{iyN} e^{-\zeta} e^{iyN} y^{-H/2}.
            (Y_{(e^{iN}.\hat F_{\infty}(s_2,\dots,s_r),W)} + \a)   
            \label{eq:twisted-limit-1}
\end{eqnarray}
By part $(c)$ of the $\SL_2$-orbit theorem (cf. equation $(4.19)$ 
in~\cite{P2}), we have
$$
   \lim_{y\to\infty} \Ad(e^{-iyN})\tilde g(y)    
   = e^{\zeta}
     \left(1+\sum_{k>0}\, \frac{(-i)^k}{k!}(\ad N)^k \tilde g_k\right)
   = e^{i\delta}
$$
Therefore, as in the proof of Theorem \eqref{theorem:limit-1} it follows
that
\beq
   \lim_{m\to\infty}\, e^{-zN}.Y_{(F(z;s_2,\dots,s_r),W)} 
   = e^{i\delta}e^{-\zeta}.Y_{(e^{iN}.\hat F_{\infty}(s_2,\dots,s_r),W)} 
   \label{eq:twisted-limit-2}
\eeq
where $\delta$ and $\zeta$ are real-analytic in $(s_2,\dots,s_r)$.
By \eqref{lemma:Deligne-2} and \eqref{corollary:Deligne-3} 
\begin{eqnarray}
      Y_{(e^{iN}.\hat F(s_2,\dots,s_r),W)}
      &=& Y(N,Y_{(\hat F(s_2,\dots,s_r),M)})            \nonumber \\
      &=& Y(N,e^{\zeta}e^{-i\delta}.Y_{(F_{\infty}(s_2,\dots,s_r),M)})
      \label{eq:twisted-limit-3} \\
      &=& e^{\zeta}e^{-i\delta}.Y(N,Y_{(F_{\infty}(s_2,\dots,s_r),M)})
         \nonumber
\end{eqnarray}
by the functoriality of Deligne's construction. Inserting 
\eqref{eq:twisted-limit-3} into \eqref{eq:twisted-limit-2} completes
the proof.
\end{proof}

\begin{remark} By virtue of the functoriality of the grading $Y(N,Y_M)$
with respect to the pair $(N,Y_M)$ and the fact that 
$Y(N,Y_M)\in\ker(\ad N)$ due to the short length of $W$, it follows that
$Y(N,Y_{(F_{\infty}(s_2,\dots,s_r),M)})$ is independent of the choice of
local coordinates.
\end{remark}

\par In connection with the proof of Theorem \eqref{theorem:main-2},
we now consider the case where $\nu$ is an admissible normal function,
on $\Delta^{*r}\subseteq\Delta^r$ with unipotent monodromy, and 
$\sing_0(\nu) =0$.  Let $(s_1(m),\dots,s_r(m))$ be a sequence
of points in $\Delta^{*r}$ which converge to $0=(0,\dots,0)$.
Let $(z_1(m),\dots,z_r(m))$ be a lifting of this sequence to
the product of upper half-planes, with the real parts of each
$z_j(m)$ restricted to an interval of finite length.  Then,
we want to compute
$$
        \lim_{m\to\infty}\, Y_{(F(z_1(m),\dots,z_r(m)),W)}
$$
where $F(z_1,\dots,z_r)$ is a lifting of the local period map to
$U^r$.  Suppose that (after passage to a subsequence)
\beq
        \lim_{m\to\infty}\, \frac{y_{j+1}(m)}{y_j(m)}\in (0,\infty)
                                                  \label{eq:condition-1}
\eeq
for $j=1,\dots,r-1$.  Then, exactly the same arguments as above show that
$$
    \lim_{m\to\infty}\, Y_{(F(z_1(m),\dots,z_r(m)),W)}
     = Y(N,Y_{(\hat F_{\infty},M)})
$$
where $N$ is any element in the monodromy cone
$
     \mathcal C = \{\, \sum_j\, a_j N_j \mid a_j>0\,\}
$.
The key point is that: 
\begin{itemize}
\item[{(a)}] By Corollary \eqref{corollary:inv-grading-split},
$
         \hat Y = Y(N,Y_{(\hat F_{\infty},M)})
$
is independent of $N$.
\item[{(b)}] Under the hypothesis of condition \eqref{eq:condition-1}, 
the element
$$
      N(y_1,\dots,y_r) = N_1 +  \frac{y_2}{y_1} N_2 
                             + \cdots + \frac{y_r}{y_1} N_r
$$
remains within a compact subset of $\mathcal C$ as $m\to\infty$.
Therefore, 
$$
     e^{iy_1 N(y_1,\dots,y_r)}.F_{\infty}
     = g(y_1)e^{iy_1 N(y_1,\dots y_r)}.\hat F_{\infty}
$$
where all the coefficients of $g$ depend real-analytically on
$N(y_1,\dots,y_r)$, since Deligne's construction \eqref{lemma:Deligne-1}
is algebraic in the pair $(N,Y_M)$.  
\end{itemize}

\par In general, by reordering the variables if necessary, one can always pass 
to some subsequence such that 
\beq
        \lim_{m\to\infty}\, \frac{y_{j+1}(m)}{y_j(m)}\in [0,\infty)
                                                  \label{eq:condition-2}
\eeq
for $j=1,\dots,r-1$.  Suppose for simplicity that 
$
        \lim_{m\to\infty}\, \frac{y_{j+1}(m)}{y_j(m)} =0 
$.
Then, the main theorem of~\cite{KNU} asserts that 
\beq
       \lim_{m\to\infty}\, Y_{(e^{iy_1 N_1 + \cdots + iy_r N_r}.F_{\infty},W)}    
                                                  \label{eq:knu-limit}
\eeq
exists (independent of any assumptions about $\sing(\nu) = 0$).

\begin{theorem}\label{theorem:limit-3} Assume that $\sing(\nu)=0$
and that $(z_1(m),\dots,z_r(m))$ is a sequence of points in $U^r$
which satisfies condition \eqref{eq:condition-2}.  Then,
$$
    \lim_{m\to\infty}\, Y_{(F(z_1(m),\dots,z_r(m)),W)}
     = Y(\sum_j\, N_j,Y_{(\hat F_{\infty},M)})
$$
\end{theorem}
\begin{proof} This is basically just the main result of \cite{KNU}
together with dependence on parameters (see the proof of the norm
estimates in \cite{KNU}) and Corollary \eqref{corollary:inv-grading-split}.
The details will appear in~\cite{BP2}.
\end{proof}

\section{Algebraicity of the Zero Locus}

\par We now prove Theorem \eqref{theorem:main-1}.  Let $\mathcal Z$
be the zero locus of an admissible normal function $\nu$ on a 
smooth complex algebraic variety $S$ which admits a smooth
compactification $\bar S$ such that $D=\bar S - S$ is a smooth
divisor.  Let $p\in D$ be an accumulation point of $\mathcal Z$,
and $(s_1,\dots,s_r)$ be local coordinates on a polydisk 
$\Delta^r\subset\bar S$ containing $p$, relative to which
$D$ is given by the equation $s_1=0$.  Let $\V\to\Delta^*\times\Delta^{r-1}$
be an admissible variation of mixed Hodge structure which
represents $\nu$ on $S\cap\Delta^r$.  Without loss of generality,
assume that $\V$ has unipotent monodromy.

\par Let $(s_1(m),\dots,s_r(m))$ be a sequence of points in $\mathcal Z$
which converge to $p$, and 
$$
           F(z;s_2,\dots,s_r):U\times\Delta^{r-1}\to\M
$$
be a lifting of the period map of $\V$, where $U$ is the upper half-plane.
Let $(z(m),s_2(m),\dots,s_r(m))$ be a lifting of $(s_1(m),\dots,s_r(m))$ to 
$U\times\Delta^{r-1}$ with the real part of $z$ restricted to an interval of 
finite length.  Then, by Theorem \eqref{theorem:limit-1}
\beq
     \lim_{m\to\infty}\, Y_{(F(z(m);s_2(m),\dots,s_r(m)),W)}
      = Y_1(0,\dots,0)                        \label{eq:limit-1}
\eeq
In particular, since the set of integral gradings is discrete,
equation \eqref{eq:limit-1} forces 
$$
          Y_{\Z} = Y_1(0,\dots,0) 
$$
to be an integral grading of $W$.  By Lemma \eqref{lemma:Deligne-2} and 
Corollary \eqref{corollary:Deligne-3}, it then follows that
\begin{itemize}
\item[{(a)}] $Y_{\Z}\in\ker(\ad N)$;
\item[{(b)}] $Y_{\Z}$ preserves the Hodge filtration 
$\hat F_{\infty}=e^{-\xi}.F_{\infty}$ of the $sl_2$-splitting 
the limit mixed Hodge structure $(F_{\infty},M)$;
\item[{(c)}] $\xi\in\ker(\ad N)\cap\lam_{(\hat F_{\infty},M)}$.
\end{itemize}

\par Let $Y_{\infty} = e^{\xi}.Y_{\Z}$.  Then, $Y_{\infty}$ preserves
$F_{\infty}$ and belongs to $\ker(\ad N)$.  Therefore, due to the short 
length of the weight filtration, there exists a unique 
$\gg_{\C}^{F_{\infty}}$-valued function $f(z;s_2,\dots,s_r)$ such that
$$
      Y_{(F(z;s_2,\dots,s_r),W)} = e^{zN}e^{\Gamma(s)}.(Y_{\infty} + f) 
$$
The local defining equation for $\mathcal Z$ near $p$ is therefore
\beq
       Y_{\Z} = e^{zN}e^{\Gamma(s)}.(Y_{\infty} + f)  \label{eq:1}
\eeq
Transposing the $e^{zN}e^{\Gamma(s)}$ factor over to the other side, we 
then obtain,
\beq
      e^{-\Gamma(s)}.Y_{\Z} = Y_{\infty} + f              \label{eq:2}
\eeq
The subalgebra ${\mathfrak q}$ is closed under the action of $\ad Y_{\infty}$.  
Consequently,
$$
        Y_{\Z} = e^{-\xi}.Y_{\infty} = Y_{\infty} + \lambda
$$
for some element $\lambda\in {\mathfrak q}$.  More properly, by equation
\eqref{eq:equivariance}, $\Lambda^{-1,-1}_{(\hat F,M)} 
= \Lambda^{-1,-1}_{(F,M)}$, wherefrom the result follows since 
$\Lambda^{-1,-1}_{(F,M)}$ is closed under $\ad Y_{\infty}$.

\par Accordingly, \eqref{eq:2} reduces to 
\beq
      e^{-\Gamma(s)}.(Y_{\infty} + \lambda) = Y_{\infty} + f   \label{eq:3}
\eeq
Again, because $Y_{\infty}$ grades $W$ and $\ad Y_{\infty}$ preserves
${\mathfrak q}$, we have
$$
       e^{-\Gamma(s)}.Y_{\infty} = Y_{\infty} + \a(s)
$$
for some holomorphic function $\a(s)$ with values in 
${\mathfrak q}\cap W_{-1}\gg_{\C}$.  Recalling that $W_{-1} \gg_{\C}$ acts 
simply transitively on the gradings of $W$, it then follows that 
equation \eqref{eq:3} simplifies to 
$$
     e^{-\Gamma(s)}.(Y_{\infty} + \lambda) = Y_{\infty}
$$
since $\gg_{\C} = {\mathfrak q}\oplus\gg_{\C}^{F_{\infty}}$ and 
$f$ takes values in $\gg_{\C}^{F_{\infty}}$.  Clearly, this 
equation is complex analytic on $\Delta^r$.  It also forces
$\lambda=0$.

\par Granting Theorem \eqref{theorem:limit-3}, the proof of 
Theorem \eqref{theorem:main-2} is identical: A sequence of points
$(s_1,\dots,s_r)$ converging to a point $p\in D$ where 
$\sing_p(\nu) = 0$ forces
$$
          Y_{\Z} = Y_{(\hat F_{\infty},M)}
$$
to be an integral grading which in the kernel of $\ad N_j$ for
each $j$.  Repeating the argument given above, it then follows
that the local defining equation for the zero locus is 
$e^{-\Gamma(s)}.Y_{\infty} = Y_{\infty}$.

\begin{remark} The above arguments also show that if $\nu$ is
an admissible normal function on $S = \bar S - D$ and $p$ is a
smooth point of $D$ such that $\sigma_{\Z,p}(\nu)$ is non-zero
torsion then $p$ can not be an accumulation point of $\mathcal Z$.
\end{remark}

%\input neron

%% input <zero-locus-smooth.bbl>
%% \bibliography{zero-locus-smooth}
%%%%% BEGIN zero-locus-smooth.bbl%%%%%%%%%%%%%%%

%%%%% END  zero-locus-smooth.bbl%%%%%%%%%%%%%%%
\end{document}